\documentclass[11pt]{preprintstyle}

\usepackage{amsfonts}
\usepackage{graphicx}

\newcommand{\curl}{\hbox{ \rm curl }}

\newcommand{\R}{\mathbb{R}}
\newcommand{\C}{\mathbb{C}}

\newcommand{\Or}{{\cal O}}

\newcommand{\setP}{{\cal P}}
\newcommand{\CalA}{{\cal A}}
\newcommand{\CalD}{{\cal D}}
\newcommand{\e}{\hbox {\rm e}}

\newcommand{\bu}{\bar{u}}
\newcommand{\bX}{{\bar{X}}}

\newcommand{\bomega}{\bar{\omega}}
\newcommand{\hv}{\hat{v}}
\newcommand{\txi}{\tilde{\xi}}
\newcommand{\tx}{\tilde{x}}
\newcommand{\td}{\tilde}
\newcommand{\pd}{{\partial}}

\newcommand{\hL}{\hat{L}}
\newcommand{\eps}{\epsilon}
\newcommand{\dps}{\displaystyle}

\newcommand\be{\begin{equation}}
\newcommand\ee{\end{equation}}
\newcommand\ba{\begin{eqnarray}}
\newcommand\ea{\end{eqnarray}}

\newcommand\nn{\nonumber}
\newcommand\bea{\begin{array}}
\newcommand\ena{\end{array}}

\begin{document}

\title{
Low-frequency stability analysis of periodic traveling-wave
solutions of viscous conservation laws in several dimensions
\footnote{Preprint: October 23, 2005}}

\date{21 October 2005} 

\author{
Myunghyun Oh 
\and Kevin Zumbrun
}

\institute{
Department of Mathematics, University of Kansas,
1460~Jayhawk~Blvd., Lawrence, KS~66047, USA 
\and Department of Mathematics, Indiana University,
Bloomington, IN~47450, USA
}

\maketitle

\begin{abstract}
We generalize work of Oh \& Zumbrun and Serre 
on spectral stability of spatially periodic traveling waves of systems
of viscous conservation laws from the one-dimensional to the multi-dimensional setting. 
Specifically, we extend to multi-dimensions the connection observed by Serre between the linearized 
dispersion relation near zero frequency of the linearized equations 
about the wave and the homogenized system obtained by slow 
modulation (WKB) approximation. 
This may be regarded as partial justification of the WKB expansion;
an immediate consequence is that hyperbolicity of the multi-dimensional homogenized system is
a necessary condition for stability of the wave. 
As pointed out by Oh \& Zumbrun in one dimension,
description of the low-frequency dispersion relation is also a first step in the determination of
time-asymptotic behavior.
\end{abstract}

\section{Introduction}\label{intro}

Nonclassical viscous conservation laws
arising in multiphase fluid and solid mechanics
exhibit a rich variety of traveling wave phenomena,
including homoclinic (pulse-type) and periodic solutions
along with the standard heteroclinic (shock, or front-type)
solutions.
Here, we investigate stability of periodic traveling waves:
specifically, the spectrum of the linearized operator about the wave.  
Our main result generalizes the works \cite{OZ.1}, \cite{Se.1}
about stability of periodic traveling waves of systems
of viscous conservation laws from the
one-dimensional to the multi-dimensional setting. 

Consider a system of conservation laws
\be
u_t + \sum_j f^j(u)_{x_j} = \sum_{j, \; k} (B^{jk}(u)u_{x_k})_{x_j},
\label{e:cons}
\ee
$u \in {\cal U} (\hbox{open}) \in \R^n, \; f^j \in \R^n, 
B^{jk} \in \R^{n \times n}, x \in \R^d, d \ge 2$, 
and a periodic traveling wave solution 
\be
u=\bar{u}(x\cdot \nu -st),
\ee
of period $X$, satisfying the traveling wave ordinary differential equation
\be
(\sum_{j,k}  \nu_j \nu_k B^{jk}(\bar u) \bar u')'= 
(\sum_j \nu_j f^j(\bar u))'-s\bar u',
\label{e:second_order}
\ee
with initial conditions
$$
\bar u(0) = \bar u(X)=:u_0.	
$$
Integrating (\ref{e:second_order}), we reduce to a first-order profile 
equation
\be
\sum_{j,k}  \nu_j \nu_k B^{jk}(\bar u) \bar u'= 
\sum_j \nu_j f^j(\bar u) -s \bar u -q
\label{e:profile}
\ee
encoding the conservative structure of the equations,
where $q$ is a constant of motion.

The one dimensional study was carried out by Oh \& Zumbrun \cite{OZ.1} 
in the ``quasi-hamiltonian'' case that the traveling-wave equation
possesses an integral of motion, and in the general case
by Serre \cite{Se.1}. 
An important contribution of Serre was to point out a 
larger connection between the linearized dispersion relation 
(the function $\lambda(\xi)$ relating spectra to wave number of the linearized
operator about the wave) near zero
and the homogenized system obtained by slow modulation approximation,
from which the various stability results of \cite{OZ.1}, \cite{Se.1}
may then be deduced.
The purpose of this paper is to extend to multiple dimensions
this important observation of Serre, 
relating the linearized dispersion relation near zero to a multi-dimensional
version of the homogenized system developed in \cite{Se.1}.
As an immediate corollary, similarly as in \cite{OZ.1}, \cite{Se.1} in the
one-dimensional case, this yields as a necessary condition for
multi-dimensional stability the hyperbolicity of the multi-dimensional
homogenized system.
In case of stability (so far not found), this relation is also the
first step in the derivation of asymptotic behavior, as in \cite{OZ.2}
in the one-dimensional case;
this we defer to a future investigation.

We here make only generic assumptions like those in \cite{Se.1},
ensuring that the set 
of periodic traveling waves is a manifold of maximal dimension
subject to the conservative properties of the equations (encoded
in form (\ref{e:profile})).
Given 
$
(a, s, \nu, q) \in {\cal U} 
\times \R \times S^{d-1} \times \R^n, 
$
(\ref{e:profile}) admits a unique local
solution $u(y; a, s, \nu, q)$ such that $u(0; a, s, \nu, q)=a$. 
Denote $X$ the period, $\omega:=1/X$ the frequency
and $M$ and $F^j$ the averages over the period:
$$
M :=\frac{1}{X} \int_0^X u(y) dy, \quad
F^j:= 
\frac{1}{X} \int_0^X 
\Big(f^j(u)- 
\sum_{k=1}^d B^{jk}(u)\omega \nu_k \partial_y u\Big) dy	
$$
when $u$ is a periodic solution of (\ref{e:profile}). Since 
these quantities are translation invariant, we consider 
the set $P$ of periodic functions $u$ that are solutions
of (\ref{e:profile}) for some triple $(s, \nu,q)$, and construct
the quotient set $\setP:=P/\cal R$ under the relation
$$
(u \; {\cal R} \; v) \Longleftrightarrow ( \exists h \in \R;\; v=u(\cdot-h)).
$$
We thus have class functions:
$$
X=X(\dot{u}),\; \omega=\Omega(\dot{u}),\;
s=S(\dot{u}),\; 
\nu=N(\dot u) , \; 
q=Q(\dot{u}),\; M=M(\dot{u}), \; F^j=F^j(\dot u),
$$
where $\dot{u}$ is the equivalence class of translates of different periodic 
functions.
Note that $\bu$ is a nonconstant periodic solution. Without loss of generality,
assume $S(\bu)=0$ and $N(\bu)=e_1$, 
%
so that (\ref{e:profile}) takes the form
$$
B^{11}(\bu)\bu'=f^1(\bu)-\bar{q}
$$
for $\bar{q}=Q(\bu)$.
Letting
$\bar{X}=X(\bu)$ and $\bar{a}=\bu(0)=u_0$, the map $(y, a, s, \nu, q) \mapsto
u(y; a,s, \nu,q)-a$ is smooth and well-defined in a neighborhood
of $(\bar{X};\bar{a},0, e_1, \bar{q})$, and it vanishes at this special point.
Here and elsewhere, $e_j$ denotes the $j$th standard Euclidean basis element.
We assume: 

(H0) $f^j$, $B^{jk}\in C^2$.

(H1) $ \hbox{Re} \; \sigma(\sum_{jk} \nu_j \nu_k B^{jk})\ge \theta>0$.

(H2) The map 
$H: \, 
\R \times {\cal U} \times \R \times S^{d-1} \times \R^n  \rightarrow \R^n$	
taking 
$(X; a, s, \nu, q)  \mapsto u(X; a, s, \nu, q)-a$
is a submersion at point $(\bar{X}; \bar{a}, 0, e_1, \bar{q})$.

As a consequence of (H0), (H2), there
is a smooth $n+d$ dimensional manifold $\setP$ of periodic solutions
$\dot u$ in the vicinity of $\bar u$, where $d$ is the spatial dimension.
On this set, one may obtain, 
rescaling by $(x,t)\to (\epsilon x, \epsilon t)$
and carrying out a formal WKB expansion as $\epsilon \to 0$ 
a closed system of $n+d$ averaged, or homogenized, equations
\ba
\label{e:wkb}
\nn
\partial_t M(\dot u) + \sum_j \partial_{x_j}(F^j(\dot u)) &=&0,  \nn \\
\partial_t (\Omega N(\dot u)) + \nabla_x  (\Omega S(\dot u))&=&0
\ea
in the $(n+d)$-dimensional unknown $\dot u$,
expected to correspond to large time-space behavior. 
For details, see Section \ref{wkb}.
The problem of stability of $\bar u$ may heuristically
be expected to be related to
the linearized equations of (\ref{e:wkb}) about the constant solution
$\dot u(x,t)\equiv u^0$, $u^0\sim \bar u$, provided that the WKB
expansion is justifiable by stability considerations.
This leads to the homogeneous degree $n+d$
linearized dispersion relation
\be\label{e:hatDelta}
\hat\Delta(\xi, \lambda):=
\det \left( 
\lambda \frac{\pd (M, \Omega N)}{\pd \dot{u}}(\dot{\bu})+
\sum_{j } i \xi_j 
\frac{\pd (F^j, 
S\Omega e_j)}{\pd \dot{u}}(\dot{\bu})  
\right )
=0.
\ee

On the other hand, one may also pursue the direct course
of linearizing PDE (\ref{e:cons}) about the stationary solution $\bar u$
and studying the spectrum of the associated linearized operator $L$.
Taking the Fourier transform in constant directions $x_j$, $j\ne 1$,
and following the general construction of \cite{G}, \cite{OZ.1},
we obtain an {\it Evans function} $D(\xi, \lambda)$, $\xi\in \R^d$,
$\lambda \in \C$, of which the zero set $(\xi, \lambda(\xi))$ determines
the linearized dispersion relation for (\ref{e:cons}), with $\lambda(\xi)$
running over the spectrum of $L$ as $\xi$ runs over $\R^d$.
For details, see Section \ref{preliminaries}. 
In particular, the low-frequency expansion of $\lambda(\xi)$ near
$(\xi, \lambda)=(0,0)$ may be expected to determine long-time asymptotic
behavior, provided that spectrum away from $\lambda=0$ has strictly
negative real part, and this in turn may be expected to derive
from the lowest order terms of the Taylor expansion of $D$.
A tedious, but fairly straightforward calculation 
following \cite{OZ.1}, \cite{Se.1} shows that
\be
\label{e:Dexpansion}
D(\xi, \lambda)= \Delta_1(\xi, \lambda) + \Or(|\xi,\lambda|^{n+2}),
\ee
where $\Delta_1$ is a homogeneous degree $n+1$ polynomial
expressed as the determinant of a rather complicated $2n\times 2n$ matrix
in $(\xi$, $\lambda)$: in particular, {\it not} in the
simple form $\det \lambda N_0+\sum_j i\xi_j N_j$ of a 
first-order hyperbolic $(n+1)\times (n+1)$ dispersion relation,
or an obvious tensor product thereof.

Our main result is the following theorem 
relating these two expansions,
generalizing the result of \cite{Se.1} in the one-dimensional case.
Define
\be\label{e:Delta}
\Delta(\xi,\lambda):=\lambda^{1-d}\hat\Delta(\xi,\lambda),
\ee
where $\hat \Delta$ is defined as in (\ref{e:hatDelta}).

\begin{theorem}\label{main}
Under assumptions (H0)--(H2),  $\Delta_1=\Gamma_0 \Delta$, i.e.,
\be
\label{e:tangent}
D(\xi, \lambda)= \Gamma_0 \Delta(\xi, \lambda)
+ \Or(|\xi, \lambda|^{n+2})
\ee
$\Gamma_0\ne0$ constant,
for $|\xi, \lambda|$ sufficiently small.
\end{theorem}

That is, up to an additional factor of $\lambda^{d-1}$,
the dispersion relation (\ref{e:hatDelta}) for the averaged system
(\ref{e:wkb}) indeed describes
the low-frequency limit of the exact linearized dispersion relation 
$$
D(\xi,\lambda)=0.
$$

The discrepancy $\lambda^{d-1}$ is an interesting and at first
glance puzzling new phenomenon in the multi-dimensional case. 
However, it is easily explained by a closer look
at the formal approximation procedure described in Section \ref{wkb}.
For, in the derivation of (\ref{e:wkb}), it was assumed that $\Omega N$
represent the gradient $\nabla_x \phi$ of a certain phase function $\phi(x,t)$.
In one dimension, this is no restriction, since we may always take
$\phi(x,t):=\int_0^x \omega (z)dz$.  However, in multidimensions,
it imposes the additional constraint
\be
\label{e:constraint}
\curl (\Omega N)\equiv 0,
\ee
which properly should be adjoined to the averaged system.

Taking the curl of the second equation of (\ref{e:wkb}),
we obtain the simple equation
$$
\partial_t \curl(\Omega N)=0,
$$
revealing at once that constraint (\ref{e:constraint}) is
compatible with the time-evolution of the system, and that
the unconstrained system possesses $(d-1)$ spurious zero
characteristices $\lambda(\xi)\equiv 0$, corresponding to
the $(d-1)$ Fourier modes in the range of the curl operator
$\hat f \to \xi \curl \hat f$, lying in $(\xi/|\xi|)^\perp$.
{\it Thus, $\Delta(\xi,\lambda)=\lambda^{1-d}\hat \Delta(\xi, \lambda)=0$ 
is exactly the linearized dispersion relation for the constrained 
averaged system (\ref{e:wkb}), (\ref{e:constraint}) 
relevant to time-asymptotic behavior.}

Theorem \ref{main} may be regarded as partial justification of the 
WKB expansion.
Roughly speaking, it states that {\it if} perturbed periodic waves
exhibit coherent behavior near the unperturbed wave $\bar u$, then
this behavior is well-described by the constrained averaged equations
(\ref{e:wkb}), (\ref{e:constraint}).
In one dimension, additional results of \cite{OZ.2} 
case give rigorous sense to this statement in the form of detailed
pointwise linear bounds under the assumption of spectral stability
of the linearized operator about the wave.
To establish a comparable {\it long time} result on behavior
in multidimensions would be a very interesting direction for future 
investigation.
Moreover, the new description of the solution given in \cite{Se.1} 
by modulation expansion might give a sufficiently good nonlinear Ansatz 
to carry out a complete nonlinear analysis, which was not done
even in the one-dimensional case.
Thus, it would be interesting to revisit also the one-dimensional setting 
of \cite{OZ.2} from this new perspective, in particular, to resolve
certain puzzling issues
in the general, non-quasi-Hamiltonian case.

An equally interesting direction for future investigation
would be to rigorously validate the WKB expansion of Section \ref{wkb}
for the closely related {\it small viscosity} problem 
\be
u_t + \sum_j f^j(u)_{x_j} = \epsilon \sum_{j, \; k} (B^{jk}(u)u_{x_k})_{x_j},
\label{e:sv}
\ee
$\epsilon \to 0$, similarly as in \cite{GMWZ.2}--\cite{GMWZ.3} for
the viscous shock case.
We note in this regard that a key ingredient in 
the \cite{GMWZ.2}--\cite{GMWZ.3} analysis,
the {conjugation lemma} of
\cite{MZ} asserting the existence of a coordinate change
converting asymptotically constant-coefficient resolvent
ODE to constant-coefficient form, {\it has a straightforward
analog in the standard Floquet construction converting
periodic- to constant-coefficient ODE}; indeed, this construction is
more natural in the periodic case.
Thus, the whole Kreiss symmetrizer construction of \cite{GMWZ.2}--\cite{GMWZ.3}
may be brought to bear also in the periodic case.
Indeed, the situation is simpler since no shock front appears;
as remarked in \cite{OZ.2}, there is a much closer analogy to
the constant-coefficient case.
Assuming successful validation of the modulation equations, 
there is also the interesting question
of secondary modulation, i.e., what happens after shock formation time
for the averaged, hyperbolic equations?

We point out further, 
as observed by Serre \cite{Se.1} in the one-dimensional
case, that the result of Theorem \ref{main} is completely 
analogous to the corresponding relation
established in \cite{ZS} for the linearized dispersion relation
associated with a perturbed viscous traveling front
$u=\bar u(x\cdot \nu -st)$, $\lim_{z\to \pm \infty}=u_\pm$,
in which the WKB expansion corresponds to matched asymptotics
joining an outer, hyperbolic solution and an inner, viscous
profile.
See \cite{Z-k}, Section 1.3, for a formal derivation in the shock
case, starting with the same rescaling $(x,t)\to (\epsilon x, \epsilon t)$;
see also the related discussion of long-time vs. small-viscosity problems in
\cite{GMWZ.1}.
Rigorous verification may be found in \cite{Z-k}, \cite{GMWZ.1}.

As an immediate consequence of Theorem \ref{main}, we obtain
the following two corollaries, yielding
a {\it necessary} condition for low-frequency
multi-dimensional spectral stability
strengthening the one-dimensional version 
obtained in \cite{OZ.1}, \cite{Se.1}.

\begin{corollary}\label{surfaces}
Assuming (H0)--(H2) and the nondegeneracy condition
\be\label{e:nondeg}
\det \left (\frac{\pd(M, \Omega N)}{\pd \dot{u}}(\dot{\bu}) \right ) \ne0,
\ee
then for $\lambda, \xi$ sufficiently small,
the zero-set of $D(\cdot, \cdot)$, corresponding to spectra of $L$, 
consists of $n+1$ characteristic surfaces:
\be
\label{e:surfaces}
\lambda_j(\xi)=-i a_j(\xi)   + \Or(\xi), \;\; j=1, \ldots, n+1,
\ee
where $a_j(\xi)$ denote the eigenvalues of 
\be\label{e:Cala}
\CalA:= \sum_{j } \xi_j \frac{\pd (F^j, S\Omega e_j)} {\pd(M, \Omega N)},
\ee
excluding $(d-1)$ identically zero eigenvalues associated
with modes not satisfying constraint (\ref{e:constraint}).
\end{corollary}

\begin{proof}
Similarly as in as in the proof of the analogous Lemma 7.5 \cite{ZS} in
the shock wave case, assuming (\ref{e:nondeg}), 
we may easily deduce (\ref{e:surfaces}) from (\ref{e:tangent}) 
using Rouch\'es Theorem.
Defining 
\be \label{cald}
\CalD^{\rho, \hat \xi}(\hat \lambda):=
\rho^{-(n+1)}D(\rho \hat \xi, \rho \hat \lambda),
\ee
for $(\rho, \hat \xi, \hat \lambda)\in \R \times S^{d-1}\times \C$,
we obtain a $d$-parameter family of analytic maps, converging as
$\rho \to 0$ to $\CalD^{0,\hat \xi}=\Delta(\hat \xi, \cdot)$.
Under assumption (\ref{e:nondeg}), 
$\CalD^{0,\hat \xi}=\Delta(\hat \xi, \cdot)\sim \hat \lambda^{n+1}$ 
as $|\lambda| \to \infty$,
hence, for $\rho$ sufficiently small, $\CalD^{\rho,\hat \xi}$ 
has $n+1$ continuously varying roots $\lambda=\hat a_j(\hat \xi, \rho)$. 
Defining $a_j(\xi):=|\xi| \hat a_j(\xi/|\xi|, |\xi|)$, we obtain the result.
\end{proof}

\begin{remark}
Evidently, $a_j(\xi)$ are smooth in $|\xi|$ for fixed $\hat \xi$,
but in general have a conical singularity at $\xi=0$ when considered as
a function of $\xi$, i.e., $\partial a_j/\partial \xi$ is discontinuous
at $\xi=0$.
\end{remark}

\begin{corollary}\label{stabcondition}
Assuming (H0)--(H2) and the nondegeneracy condition (\ref{e:nondeg}),
a necessary condition for low-frequency spectral stability of
$\bar u$, defined as
Re $\lambda \le 0$ for $D(\xi,\lambda)=0$, $\xi\in \R^d$, and
$|\xi,\lambda|$ sufficiently small, 
is that the averaged system (\ref{e:wkb}) be ``weakly hyperbolic''
in the sense that it possesses a full set of real characteristics
$\hat \lambda_j(\xi)$ for each $\xi\in \R^d$, i.e.,
the eigenvalues of 
$$
\CalA = \sum_{j } \xi_j \frac{\pd (F^j, S\Omega e_j)} {\pd(M, \Omega N)}
$$
are real.
\end{corollary}

\begin{remark}\label{onedim}
Condition (\ref{e:nondeg}), or equivalently
$(\partial/\partial \lambda)^{n+1}D(0,0)\ne 0$,
is a necessary condition for one-dimensional 
linearized stability \cite{OZ.2}, while hyperbolicity is
necessary for stability of the homogenized system
linearized about a constant state.  Thus, Corollaries
\ref{surfaces} and \ref{stabcondition} are 
analogous to results of \cite{ZS} in the shock wave case,
stating that, given one-dimensional stability, 
stability of the inviscid equations linearized about an ideal
shock is necessary for multi-dimensional stability of a viscous shock wave.
\end{remark}

Finally, we mention that, though the averaged system may in some
cases be hyperbolic \cite{OZ.1}, so far, 
only unstable periodic traveling-wave
solutions have been found for viscous conservation laws.  
However, essentially only the single, $2 \times 2$ model of van der 
Waals gas dynamics with viscosity--capillarity in one dimension 
has so far been considered in detail \cite{OZ.1}, \cite{Se.1}--\cite{Se.2}, 
and we see no obvious reason why a stable wave should not 
exist for other models.  It would be extremely interesting to 
either find such an example, with the associated
rich behavior described by the modulation equations, or
show that it can in no case exist.
As suggested by Serre \cite{Se.3}, a useful starting point
might be to consider whether the averaged system (\ref{e:wkb})
might ever possess an entropy.

{\bf Plan of the paper}  In Section \ref{preliminaries}, we
recall the basic Evans function construction of \cite{G}.
In Section \ref{evans}, we carry out the expansion (\ref{e:Dexpansion}),
and in Section \ref{wkb}, the multi-dimensional WKB expansion (\ref{e:wkb}).
Finally, in Section \ref{proof}, we carry out the proof of
Theorem \ref{main} by a calculation similar to the one used by 
Serre \cite{Se.1} to treat the one-dimensional case.

\section{Preliminaries}\label{preliminaries}

Without loss of generality taking $S(\bu)=0$, $N(\bu)=e_1$,
$\bar u=\bar{u}(x_1)$ represents a stationary solution. 
Linearizing (\ref{e:cons}) about 
$\bar{u}(\cdot)$, we obtain
\be
v_t = Lv := \sum(B^{jk}v_{x_k})_{x_j}-\sum(A^j v)_{x_j},
\label{e:lin}
\ee
where coefficients
\be
B^{jk} := B^{jk}(\bu), \;\;\;
A^jv:= Df^j(\bu)v-(DB^{j1}(\bu)v)\bu_{x_1}
\label{e:coeff}
\ee
are now periodic functions of $x_1$. 

Taking the Fourier transform in the transverse coordinate $\td{x}=
(x_2, \cdots, x_d)$, we obtain
\ba
\hv_t = L_{\td{\xi}}\hv 
& = &(B^{11} \hv_{x_1})_{x_1} -(A^1 \hv)_{x_1}
	+i (\sum_{ j\ne 1} B^{j1} \xi_j) \hv_{x_1} 	\nn \\ 
& + & i(\sum_{k\ne 1}B^{1k} \xi_k \hv)_{x_1} 
	- i \sum_{j\ne 1}A^j \xi_j \hv 
	- \sum_{j\ne 1, k\ne 1} B^{jk}\xi_k \xi_j \hv,
\label{e:fourier}
\ea
where $\td{\xi}=(\xi_2, \cdots, \xi_d)$ is the transverse frequency
vector.
The Laplace transform in time $t$ leads us to study the family of
eigenvalue equations
\ba
0=(L_{\txi}- \lambda) w
&=& (B^{11}w')'-(A^1 w)'+i \sum_{ j\ne 1}B^{j1} \xi_j w'
	+i(\sum_{k \ne 1} B^{1k} \xi_k w)'	\nn \\
&-& i\sum_{j \ne 1}A^j \xi_j w - \sum_{j\ne 1, k\ne 1}B^{jk}
	\xi_k \xi_j w - \lambda w,
\label{e:laplace}
\ea
associated with operators $L_{\txi}$ and frequency $\lambda \in \C$,
where `$'$' denotes $\partial/\partial x_1$.
Clearly, a necessary condition for stability of (\ref{e:cons}) is
that (\ref{e:laplace}) have no $L^2$ solutions $w$ 
for $\txi \in \R^{d-1}$ and Re $\lambda >0$. 
For solutions of
(\ref{e:laplace}) correspond to normal modes $\hv(x,t)
=\e^{\lambda t} \e^{i \txi \cdot \tx} w(x_1)$ of (\ref{e:lin}).

The difficulty of our problem is due to accumulation at the
origin of the essential spectrum of the linearized operator $L$
about the wave as in the one dimensional case. 
Multidimensional stability concerns
the behavior of the perturbation of the top eigenvalue, 
$\lambda=0$ under small perturbations in $\txi$.
For study this stability, we use Floquet's theory and an Evans function 
\cite{G} which
not only depends on $\lambda$ but also on $\xi_1$ which corresponds
to the phase shift and $\txi$.
To define the Evans function, we choose a basis 
$\{w^1(x_1, \txi, \lambda), \ldots,
w^{2n}(x_1, \txi, \lambda) \}$ of the kernel of $L_{\txi}-\lambda$, which is
analytic in $(\txi, \lambda)$ and is real when $\lambda$ is real, for details
see \cite{OZ.1,Se.1}. 
Now we can define the Evans function by
\be
D(\lambda, \xi_1, \txi) := \left | 
\begin{array}{c}
  w^l(X, \txi, \lambda)-\e^{iX\xi_1}w^l(0, \txi, \lambda)  \\
  (w^l) '(X, \txi, \lambda)-\e^{iX\xi_1}(w^l) '(0, \txi, \lambda) \\
\end{array}
\right |_{1 \le l \le 2n}
\label{e:evans}
\ee
where $\xi_1 \in \R$. Note that $X\xi_1$ is exactly $\theta$ in \cite{Se.1}. 
We remark that $D$ is analytic everywhere, with associated
analytic eigenfunction $w^l$ for $1 \le l \le 2n$. 
A point $\lambda$ is in the spectrum
of $L_{\txi}$ if and only if  $D(\lambda, \xi)=0$ with $\xi=(\xi_1, \txi)$.

\begin{example}
In the constant-coefficient case
\ba
&& B^{11}w''-A^1 w'+i \sum_{ j\ne 1}B^{j1} \xi_j w'
	+i\sum_{k \ne 1} B^{1k} \xi_k w' 	\nn \\
&-& i\sum_{j \ne 1}A^j \xi_j w - \sum_{j\ne 1, k\ne 1}B^{jk}
	\xi_k \xi_j w - \lambda w=0,
\label{e:const}
\ea
an elementary computation yields
$$
D(\lambda, \xi)= \Pi_{l=1}^{2n} (\e^{\mu_l(\lambda, \txi)T}-\e^{i \xi_1 T})
$$
where $\mu_l, l=1, \ldots, 2n$, denote the roots of the characteristic equation
\ba
&& ( \mu^2 B^{11}+ \mu(-A^1 +i \sum_{ j\ne 1}B^{j1} \xi_j 
	+i\sum_{k \ne 1} B^{1k} \xi_k ) 	\nn \\
&-& ( i\sum_{j \ne 1}A^j \xi_j  + \sum_{j\ne 1, k\ne 1}B^{jk}
	\xi_k \xi_j  + \lambda I )) \bar{w} =0,
\label{e:const_char}
\ea
where $w=\e^{\mu x_1} \bar{w}$. The zero set of $D$ consists of all 
$\lambda$ and $\xi_1$
such that 
$$
\mu_l(\lambda, \txi)=i \xi_1 (\hbox{mod} 2 \pi i /X)
$$
for some $l$. 
Setting $\mu= i \xi_1$ in (\ref{e:const_char}), 
we obtain the dispersion relation
\be
\hbox{det} (- B^{\xi} - i A^{\xi} - \lambda I)=0
\label{e:const_dispersion}
\ee
where $A^{\xi}= \sum_j A^j \xi_j$ and $B^{\xi}= \sum_{j, k} B^{jk} \xi_k \xi_j$.
\end{example}

\begin{remark}
If $(\lambda, \txi)=(0,0)$ then (\ref{e:const_char}) reduces to 
$$
\mu( (B^{11})^{-1}A^1- \mu) \bar{w} =0
$$
giving $n$ nonzero roots
$$
\mu=s_j, \;\; \bar w =t_j,
$$
where $s_j, t_j$ are eigenvalues and eigenvectors for the matrix 
$(B^{11})^{-1}A^1$, and an $n$-fold root $\mu=0$. 
Thus, $D(0,0)=0$ in the above example. 
We shall see later that this holds also in the general 
variable-coefficient case. 
\end{remark}

\begin{remark}
In the constant coefficient case, (\ref{e:const_dispersion}) yields expansions
\be
\lambda_j(\xi)=0- i a_j(\xi)  + \Or(|\xi|), \;\;j=1, \ldots, n,
\label{e:const_lambda}
\ee
for the $n$ roots bifurcating from $\lambda(0)=0$, where $a_j$ denote the
eigenvalues of $A^{\xi}$. Thus we obtain the necessary stability condition of 
hyperbolicity, $\sigma(A^{\xi})$ real. 
\end{remark}

\section{Evans function calculations}\label{evans}

Motivated by the example, we now find linearized dispersion relations
for the variable-coefficient Evans function in the low-frequency limit.
From now on, coordinatize $\nu$ in the vicinity of $e_1$ by
\be\label{e:delta}
\nu=:\frac{(1, \delta^2, \ldots, \delta^d)}{\sqrt{1+ |\delta|^2}},
\ee
$\delta=(\delta^2, \ldots, \delta^d)\in \R^{d-1}$.
Note that differentiation of (\ref{e:delta}) yields  
$\partial \nu= (0, \delta)$. 

\begin{lemma}[\cite{Se.1}]
Assumption (H2) is equivalent to 

(H2)$'$
$$
\R^n=\hbox{span} \{ \left[ \frac{\pd u}{\pd s} \right ], 
\left[ \frac{\pd u}{\pd \delta} \right ], 
[w^2], \ldots, [w^{2n}], \bu'(0) \}.
$$
\end{lemma}

\begin{proof}
Immediate, using $[w_1]=0$; see \cite{Se.1} for the one-dimensional case.
\end{proof}

\subsection{Variational relations}

If $\txi=0$, (\ref{e:laplace})  equation becomes 
\be
(L_0-\lambda)w=(B^{11}w')'-(A^1 w)' - \lambda w=0,
\label{e:simple}
\ee
which is associated with the one
dimensional stability problem studied in \cite{OZ.1,Se.1}. 
Recall that $\bu$ is $\bX$-periodic in $x_1$ and 
the functions $w^1(x_1, \txi, \lambda), \ldots,
w^{2n}(x_1, \txi, \lambda)$ are in the basis of
the kernel of $L_{\txi}-\lambda$. 
Following \cite{OZ.1}, we normalize 
\ba
\label{normal}
\nn
w^j(0,\tilde \xi, \lambda)&=&e_j, 
\quad (w^j)'(0,\tilde \xi, \lambda)=(B^{11})^{-1}A^1 e_j; \\
\nn
\quad
w^{n+j}(0,\tilde \xi, \lambda)&=&0, \quad (w^{n+j})'(0,\tilde \xi, \lambda)=-(B^{11})^{-1}e_j
\ea
for $j=1, \dots, n$ and all $(\tilde \xi, \lambda)$,
giving in particular
\be \label{e:zerovalues}
\hat L w^j(0,\tilde \xi, \lambda)=0, \,\;\;
\hat L w^{n+j}(0,\tilde \xi, \lambda)=e_j
\;\;\hbox{ for} \,\;\;
j=1, \dots, n.
\ee

Plug the Taylor expansion of $w(x_1, \txi, \lambda)$ at the origin of $(\lambda, \txi)$
\ba
w^l(\cdot, \txi, \lambda)&=&w^l(\cdot, 0,0)+\lambda w^l_\lambda(\cdot, 0,0) 
+ \sum_{j\ne1}w^l_{\xi_j}(\cdot, 0,0) \xi_j \\		\nn
&& + \frac{1}{2}( \lambda^2 w^l_{\lambda \lambda}
(\cdot,0,0)+ 2\lambda \sum_{j\ne 1} w^l_{\lambda \xi_j}(\cdot,0,0) \xi_j 
+ \sum_{j \ne 1, k \ne 1} w^l_{\xi_j \xi_k}(\cdot,0,0) \xi_k \xi_j) + \ldots
\ea
into (\ref{e:laplace}) to find the identies:
\be
(\hat{L}w^l)'=0, \quad (\hat{L} w^l_{\lambda})'=w^l, 
\quad (\hat{L} w^l_{\lambda\lambda})' =2 w^l_{\lambda}
\label{e:identity}
\ee
and
\ba
\nn
(\hat{L} w^l_{\lambda\xi_j})' &=& (iA^j w^l_\lambda 
-iB^{j1}(w_{\lambda}^l)'-i(B^{1j}w_\lambda^l)' 
+ w^l_{\xi_j}), \\\nn
(\hat{L} w^l_{\xi_j\xi_k})' &=&
(iA^j w^l_{\xi_k}+ iA^k w^l_{\xi_j}
-iB^{j1}(w_{\xi_k}^l)'-i(B^{1j}w_{\xi_k}^l)' 
-iB^{k1}(w_{\xi_j}^l)'-i(B^{1k}w_{\xi_j}^l)' 
+2B^{jk}w^l
), 
\label{e:mixedidentity}
\ea
where $\hat{L}w = B^{11} w' -A^1 w$, and also: 
\be
(\hL w^l_{\xi_j})' = i ( A^jw^l-B^{j1}(w^l)' - (B^{1j}w^l)' ),
\quad j\ne 1
\label{e:identity_1}
\ee
and 
\be
(\hat{L}w^1_{\xi_j})'=i (f^j(\bu)-B^{j1}(\bu)\bu'-B^{1j}(\bu)\bu')',
\quad j\ne 1
\label{e:extra_1}
\ee
by using the definition of $A^j$ in (\ref{e:coeff}).
{In the Laplacian case $B^{jk}=\delta^j_k$}, the latter identity
simplifies to
\be
(\hat{L}w^1_{\xi_j})'=i f^j(\bu)'.
\label{e:Lextra_1}
\ee

Note that the functions $w^1(x_1, 0, 0), \ldots,
w^{2n}(x_1, 0, 0)$ are in the basis of
the kernel of $L_0$.
We omit $(\cdot, 0,0)$ hereabove and denote  $[w]=w(X)-w(0)$.
We also have $w^1=\bu'$, $\hL w^1=0$, 
$[w^1]=0$, $\int_0^{\bX} w^1 dx_1=0$, and  
moreover, 
\be
[\hat{L}w^l]=0, \quad [\hat{L} w^l_{\lambda}]=\int_0^{\bX} w^l dx_1, 
\label{e:extra}
\ee
and 
\be
[\hat{L} w^1_{\lambda}]=0, \quad
[\hat{L} w^1_{\xi_j}]=0, \quad
\ee
\ba
\nn
\label{secondorder}
[\hat{L} w^1_{\lambda\lambda}]& = & 2 \int_0^{\bX} w^1_{\lambda} dx_1,  \\
\nn
[\hat{L} w^1_{\lambda \xi_j}]& = & 	
\int_0^{\bX} 
\Big(iA^j w^1_\lambda -iB^{j1}(w_{\lambda}^l)'
-i(B^{1j}w_\lambda^1)' + w^1_{\xi_j}\Big)
 dx_1,   
\\
\nn
[\hat{L} w^1_{\xi_j \xi_k}]& = & 	
\int_0^{\bX } 
\Big(iA^j w^1_{\xi_k}+ iA^k w^1_{\xi_j}
-iB^{j1}(w_{\xi_k}^1)'-i(B^{1j}w_{\xi_k}^1)' \\
\nn
& \quad &-iB^{k1}(w_{\xi_j}^1)'-i(B^{1k}w_{\xi_j}^1)' +2B^{jk}w^1\Big) 
dx_1 .  \\
\ea

{In the Laplacian case $B^{jk}=\delta^j_k$}, the last two identities
simplify considerably, to
\ba
\label{Lmixed}
\nn
[\hat{L} w^1_{\lambda \xi_j}] &=&
\int_0^{\bX} (iA^j w^1_\lambda+ w^1_{\xi_j})dx_1,\\
\nn   
[\hat{L} w^1_{\xi_j \xi_k}] &=&  	
\int_0^{\bX } 
\Big(iA^j w^1_{\xi_k}+ iA^k w^1_{\xi_j} \Big)
dx_1 .  \\
\ea

\subsection{Connection to traveling-wave variations}
From (\ref{normal}), we find easily that 
\be
\label{init}
w^j(\cdot, 0,0)=\partial u/\partial a_j|_{\bar u},
\qquad
w^{n+j}(\cdot, 0,0)=\partial u/\partial q_j|_{\bar u} \;\;  \hbox{for}\;\; j=1, \ldots, n. 
\ee
For example, taking the variation of traveling wave equation
(\ref{e:profile})
with respect to $q_j$, we
find that $z=\pd u/\pd q_j$ satisfies  
$$B^{11} z'-A^1z= -e_j$$
with $z(0)=0$, so that $L_0z=0$ and $z'(0)= -(B^{11})^{-1}e_j$ as claimed.   

Further (see \cite{OZ.1}, \cite{Se.1}), 
$$
w^1_{\lambda}= -\pd u / \pd s+ \sum_{n+1}^{2n} \alpha^l w^l
$$ 
for $\alpha \in \R^{2n}$, since $L_0 (-\partial  u/\partial s)
=L_0 w^1_\lambda(\cdot, 0,0)= \bar u'$ 
and $w^1_\lambda(0,0,0)=(\partial u/\partial s)(0)=0$,
and, similarly, for $j \ne 1$, 
using $L_0w^1_{\xi_j}(\cdot, 0,0)= L_0(i\partial u/\partial \delta^j)(\cdot, 0,0)=
if^j(\bar u)'$ and $w^1_{\xi_j}(0,0,0)=\partial u/\partial \delta^j(0)=0$,
$$
w^1_{\xi_j}(\cdot, 0,0)= i \pd u / \pd \delta^j
+ \sum_{n+1}^{2n} \beta^l_j w^l.
$$ 

Alternatively, 
\be \tilde w^1_\lambda(\cdot,0,0)= -\pd u / \pd s, \qquad
\tilde w^1_{\xi_j}(\cdot,0,0)= i \pd u / \pd \delta^j
\label{e:var_rels}
\ee
for 
\be
\tilde w^1:= w^1-\lambda \sum_{n+1}^{2n} \alpha^l w^l
- \sum_{\ell= n+1}^{2n} \sum_j \xi_j \beta^l_j w^l,
\ee
with $\tilde w^1(0,0)$ still equal to $\bar u'$.
We hereafter substitute $\tilde w^1$ for $w^1$ everywhere it
appears, as we are free to do.  (Recall, $w^\ell$ can be
an arbitrary basis of the kernel of $L$.)

\subsection{Reduction of the leading part}
We rewrite the Evans function (\ref{e:evans}) as
\be
D(\lambda, \xi_1, \txi) := \left |
\begin{array}{c}
 [ w^l(\txi, \lambda)]+ (1-\e^{i\bX\xi_1}) w^l(0, \txi, \lambda)  \\
\left [ (w^l) ' (\txi, \lambda) \right ]  +(1-\e^{i\bX\xi_1})(w^l) '(0, \txi, \lambda) \\
\end{array}
\right |_{1 \le l \le 2n}
\label{e:evans_re}
\ee
and then multiply the second row in (\ref{e:evans_re}) by $B^{11}$
and then subtract $A^1$ times the first one
\be
( \det B^{11})  D(\lambda, \xi_1, \txi) := \left |
\begin{array}{c}
 [ w^l(\txi, \lambda)]+ (1-\e^{i\bX\xi_1}) w^l(0, \txi, \lambda)  \\
\left [  \hat{L} w^l (\txi, \lambda) \right ] +(1-\e^{i\bX\xi_1})(\hat{L} w^l) (0, \txi, \lambda) \\
\end{array}
\right |_{1 \le l \le 2n}.\\
\label{e:evans_final}
\ee

{\it At this point, we restrict for readibility to the Laplacian case
$B^{jk}=\delta^j_k$}.  The general case goes similarly.
Then the Evans function $D(\lambda, \xi_1,0)$ becomes
\ba \nn
&& \Gamma_0 \det 
\left(
\begin{array}{ccc}
c(\xi,\lambda) & [w^2]  & \ldots  \\
C(\xi,\lambda) & \;\;
C^2(\xi, \lambda) & \ldots
\end{array}
\right)	\\ 
&+& \Or(|\lambda|^{n+2}+|\xi_1|^{n+2})
\label{e:evans_variable}
\ea
with a nonzero number $\Gamma_0$ (for details, see \cite{OZ.1,Se.1}),
where
\ba
\nn
c(\xi,\lambda)&=&
\lambda[w^1_{\lambda}] +\sum_{j \ne 1}\xi_j[w^1_{\xi_j}]- i \bX \xi_1 \bu'(0) \\
&=&
-\lambda[\pd u/\pd s] +\sum_{j \ne 1} i\xi_j[\partial u/\partial \delta_j] 
-i \bX \xi_1 \bu'(0) \\
\ea
is a homogeneous degree one polynomial,
\ba
\nn
C(\xi,\lambda)&:=&	
\frac{1}{2} \lambda^2 [\hat L w^1_{\lambda \lambda}]
+ \lambda \sum_{j\ne 1} \xi_j [\hat L w^1_{\lambda \xi_j}]
+\frac{1}{2} \sum_{j,k\ne 1} \xi_j \xi_k [\hat L w^1_{\xi_j \xi_k}]
\nn \\
& & \quad 
-i \bX \xi_1 \lambda (\hat{L} w^1_\lambda ) (0) 
-
i \bX \xi_1 \sum_{j\ne 1} \xi_j (\hat{L} w^1_{\xi_j} ) (0) 
\nn \\
&=&
-\lambda^2 \int_0^{\bX} (\pd u/\pd s) dx_1
+i\lambda \sum_{j\ne 1} \xi_j
\int_0^{\bX} (A^j (-\pd u/\pd s)+ (\pd u/\pd \delta^j)) dx_1
\nn \\
& & \quad
- \sum_{j,k\ne 1} \xi_j \xi_k
\int_0^{\bX } 
A^j (\pd u /\pd \delta^k) 
dx_1 
\nn \\
& & \quad 
+i \bX \xi_1 \lambda (\hat{L} \pd u/ \pd s) (0) 
+ \bX \xi_1 \sum_{j\ne 1} \xi_j (\hat L \pd u/\pd \delta^j)(0)
\nn \\
&=&
-\lambda^2 \int_0^{\bX} (\pd u/\pd s) dx_1 
+i\lambda \sum_{j\ne 1} \xi_j
\Big(-(\pd/\pd s)\int_0^{\bX} f^j (u) dx_1 + 
(\pd /\pd \delta^j)\int_0^{\bX} u dx_1
\Big)
\nn \\
& & \quad
-\sum_{j,k\ne 1} \xi_j \xi_k
(\pd/\pd \delta^k)\int_0^{\bX } 
f^j (u) dx_1 
\nn \\
& & \quad 
+i \bX \xi_1 \lambda (\hat{L} \pd u/ \pd s) (0) 
+ \bX \xi_1 \sum_{j\ne 1} \xi_j (\hat L \pd u/\pd \delta^j)(0)
\ea
is a homogeneous degree two polynomial, and
\be
C^\ell(\xi,\lambda):=
\lambda [\hL w^l_{\lambda}]
+\sum_{j\ne1} \xi_j [\hat L w^l_{\xi_j}]
- i \bX \xi_1\hat{L} w^l(0) 
\ee
are homogeneous degree one polynomials given by
\be
(\pd/\pd a_j)
\Big( \lambda \int_0^{\bX} u(x_1) dx_1
+\sum_{j\ne1} \xi_j \int_0^\bX i f^j(u(x_1)) dx_1
- i \bX \xi_1(u'(0)-f^1(u(0)) \Big)
\ee
for $\ell=j=2, \ldots, n$, and
\be
(\pd/\pd q_j)
\Big( \lambda \int_0^{\bX} u(x_1) dx_1
+\sum_{j\ne1} \xi_j \int_0^\bX i f^j(u(x_1)) dx_1
- i \bX \xi_1(u'(0)-f^1(u(0)) \Big)
\ee
for $\ell=n+j=n+1, \ldots, 2n$.

Thus, the leading order part of $D$ near $(\xi,\lambda)=(0,0)$
is the homogeneous degree $(n+1)$ polynomial
\ba \nn
\Delta_1(\xi,\lambda) &:=& \Gamma_0 \det 
\left(
\begin{array}{ccc}
c(\xi,\lambda) & \ldots  \\
C(\xi,\lambda) & \;\;
C^2(\xi, \lambda) & \ldots
\end{array}
\right)	\\ 
\label{e:leading}
\ea
with $c$, $C$, $C^\ell$ defined as above.
In particular, the Evans function has a zero of order
$n+1$ at $(\lambda, \xi)=(0,0)$.

\section{Slow modulation approximation}\label{wkb}
Next, we carry out a multi-dimensional version of the slow modulation (WKB) 
expansion in \cite{Se.1}.
Rescale $(x,t) \mapsto (\eps x, \eps t)$ in (\ref{e:cons}) to obtain
\be
u_t + \sum_j f^j(u)_{x_j} = \eps \sum_{j, \; k} (B^{jk}(u)u_{x_k})_{x_j}.
\label{e:rescaled_cons}
\ee
Let 
\be
u^{\eps}(x,t)=u^0 \left (x,t, \frac{\phi(x,t)}{\eps} \right )
+ \eps u^1 \left (x,t, \frac{\phi(x,t)}{\eps} \right ) + \cdots,
\label{e:homo_sol}
\ee
where $y \mapsto u^0(x,t,y)$ is a periodic function with $\partial_{x} \phi \ne 0$.
We plug (\ref{e:homo_sol}) into (\ref{e:rescaled_cons}) and consider the 
equations obtained by equating coefficients at successive powers of $\epsilon$.


At order $\eps^{-1}$, we have
\ba	\nn
-s \partial_y u^0 &+& \sum_j \omega \nu_j \partial_y (f^j(u^0))
- \partial_y(\sum_{j,k} \omega^2 \nu_j \nu_k B^{jk}(u^0) \partial_y u^0)= 0,
\ea
with 
\be 
\dps s:= - \frac{\partial_t \phi}{|\partial_{x} \phi|},
\quad
\nu:= \frac{\partial_x \phi}{|\partial_{x} \phi|},
\quad
\omega:= |\partial_{x} \phi|,
\ee
which may be recognized as the traveling-profile equation after
rescaling $y\to \omega y$.  That is, $u^0(y)=\bar u(\omega y)$
for a periodic profile of period $X=\omega^{-1}$, hence
$u^0$ is periodic of {\it period one}, as described in \cite{Se.1}.
The quantities $\omega(x,t)$, $s(x,t)$, $\nu(x,t)$ are the local frequency,
speed, and direction of the modulated wave.

At order $\eps^0$, we have
\ba	\nn
\partial_t u^0 + 
\sum_{j=1}^d \partial_{x_j}\Big(f^j(u^0)- 
\sum_{k=1}^d B^{jk}(u^0)\omega \nu_k \partial_y u^0\Big)	
&=&\partial_y (\ldots).
\ea
Taking the average with repect to $y$, 
and rescaling with $y:= \omega y$, we obtain 
\be
\partial_t M(u^0)+ \sum_{j = 1} \partial_{x_j}  F^j (u^0)
=0
\label{e:homo1}
\ee
where 
\be
\label{Fj}
F^j (u^0)=
\frac{1}{X} \int_0^X 
\Big(f^j(u^0)- 
\sum_{k=1}^d B^{jk}(u^0)\omega \nu_k \partial_y u^0\Big) dy	
\ee
is the averaged flux along orbit $u^0$ (now rescaled to actual
period $\bX$), with
$$
\sum_j \nu_j F^j= (SM+Q)(u^0),
$$
by the profile equation.
{In the Laplacian case $B^{jk}=\delta^j_k$}, (\ref{Fj})
simplifies to
\be
\label{LFj}
F^j (u^0)=
\frac{1}{X} \int_0^X 
(f^j(u^0)-\nu_j (u^0)\partial_y u^0) dy.
\ee
We have an additional $d$ equations
\be
\partial_t (\Omega N) (u^0)+\partial_{x}(\Omega S(u^0))=0
\label{e:homo2}
\ee
from the Schwarz identity 
$\partial_t \partial_x \phi=\partial_x\partial_t \phi$,
where $d$ is the dimension of the spatial variable $x$.
(Note: $(\Omega, N)$ may be regarded as polar
coordinates for $\Omega N$.)

Combining, we obtain finally the closed homogenized system 
\be
\label{e:homo}
\partial_t (M, \Omega N) + \sum_j \partial_{x_j}(F^j, \Omega Se_j)
\ee
of the introduction, consisting of $n+d$ equations in $n+d$ unknowns. 
As discussed in the introduction, this should be supplemented with
the constraint
\be
\curl (\Omega N)\equiv 0
\ee
coming from the relation $\Omega N= \nabla_x \phi$.

\section{Proof of the main theorem}\label{proof}
We now carry out the proof of Theorem \ref{main},
{restricting for readibility to the Laplacian case 
$B^{jk}=\delta^j_k$}.  The general case follows similarly.
We want to see that 
the leading order part $\Delta_1$ of $D$, defined
in (\ref{e:leading}), is given by a (nonzero) constant multiple of 
$\lambda^{1-d}$ times
\be
\hat \Delta(\xi,\lambda)=
\det \left( 
\lambda \frac{\pd (M, \Omega N)}{\pd \dot{u}}(\dot{\bu})+
\sum_{j } i \bar{\omega} \bX \xi_j 
\frac{\pd (F^j, 
S\Omega e_j)}{\pd \dot{u}}(\dot{\bu})  
\right ),
\label{e:delta_multid}
\ee
where $\dot {\bar u}$ denotes the orbit class of $\bar u$,
with $\bar\nu=(1,0,\dots,0)= N(\dot{\bar u})$.

Recall the assumption (H2)$'$ and (\ref{e:delta}). The tangent
space to $P$ at $\bu$ is the $\beta$-projection of 
the kernel of 
$$
Z(\beta, \delta, \gamma):=\beta^0 \left [\frac{\partial u}{\partial s} \right ] 
+ \sum_1^{2n} \beta^l [w^l] + 
\sum_{j=2}^d \delta^j 
\left [\frac{\partial u}{\partial \delta^j} \right ] 
+\gamma \bu'(0),
$$
and the tangent space to $\cal P$ at $\dot{\bu}$ is the $\beta$-projection of 
the kernel of 
$$
{\cal Z}(\beta^0, \beta^2, \ldots, \beta^{2n}, 
\delta^2, \ldots, \delta^{d}, \gamma):=
\beta^0 \left [\frac{\partial u}{\partial s} \right ] 
+ \sum_2^{2n} \beta^l [w^l]
+\sum_{j=2}^d \delta^j 
\left [\frac{\partial u}{\partial \delta^j} \right ] 
+\gamma \bu'(0).
$$
We relabel $\beta^1$ for $\beta^0$ since we will not use $w^1$ hereafter. 
Thus, 
\be \label{e:cal_Z}
{\cal Z}(\beta, \delta, \gamma):= \beta^1 \left [
\frac{\partial u}{\partial s} \right ] + \sum_2^{2n} \beta^l [w^l]
+\sum_{j=2}^d \delta^j 
\left [\frac{\partial u}{\partial \delta^j} \right ] 
+\gamma \bu'(0), \;\; \beta \in \C^{2n}, \; \delta\in \C^{d-1} \; \gamma \in \C.
\ee

We easily compute (see (\ref{e:zerovalues}) for $dQ$) the differentials
\ba
&&dX \cdot (\beta, \delta, \gamma)=\gamma=\partial X, \;\; 
dS \cdot (\beta, \delta, \gamma)=\beta^1=\partial S, \;\;
dN\cdot (\beta, \delta, \gamma)=
(0, \delta)^T=\partial N,
\\	\nn
&& dQ \cdot (\beta, \delta, \gamma)=
(\beta^{n+1}, \ldots, \beta^{2n})^T
=-\sum_2^{2n} \beta^l \hL w^l=\partial Q,  \;\;
d\Omega \cdot (\beta, \delta, \gamma)=- \bar{\omega}^2 \gamma=\partial \Omega,
\label{e:diff1}
\ea
where, following \cite{Se.1}, we use the notation $\partial G$ to indicate
the extension to $\C^{2n+d}$ of a differential $dG$ defined on the 
kernel of ${\cal Z}$.  
(Note: this includes the extension from real to complex values,
of which we shall later make important use
in parametrizations (\ref{e:comp1}) and (\ref{e:comp2}).)

Likewise, $XM=\int_0^X u(y) dy$ gives
\be
d(XM)\cdot (\beta, \delta, \gamma)=
\gamma\bu(0)+\beta^1\int_0^{\bar{X}}
\frac{\partial u}{\partial s} dy+
\sum_2^{2n}\beta^l \int_0^{\bar{X}}w^l dy +
\sum_2^d \delta^j
\int_0^{\bar{X}} \frac{\partial u}{\partial \delta^j} dy
 =\partial(XM),
\ee
thus determining
\be\label{e:Mvar}
\partial M = \bar \omega (\partial (X M)-{\bar M} \partial X).
\ee

Finally, variations
\be
\label{Fjvar}
d(XF^j)(u)  =
d \int_0^X 
(f^j(u)-\nu_j (u)\partial_y u) dy= \partial (XF^j) (u)
\ee
for $j\ne1$ have the simple form
\be
\label{XFjvarsimple}
\partial(XF^j) (u)=
\pd X f^j(u)(0) +
\int_0^X \partial f^j(u) dy,
\ee
by (\ref{LFj}) and periodicity of $u$, hence
\ba
\label{Fjvarsimple}
\partial F^j (u)&=& \omega(\partial( XF^j)- (\partial X)F^j(u))
\nn \\
&=&
\omega \gamma  f^j(u)(0) +
\omega \int_0^X \partial f^j(u) dy 
- \omega \gamma F^j(u).
\ea

To find the variation for $F^1$, note that,
by the first-order traveling wave system (\ref{e:profile}), 
$$
\sum_j \nu_j F^j = MS + Q,
$$
so that $\sum_j (\partial \nu_j)F^j + \sum_j \nu_j (\partial F^j)=
\partial (MS+Q)$, hence, for $\nu=(1,0,\dots,0)$, $s=0$,
\be
\label{F1varsimple}
\partial F^1 (u)
= \partial (M S+ Q) - \sum_j \delta^j F^j(u)
= M \partial S+ \partial Q - \sum_{j\ne 1} \delta^j F^j(u).
\ee

We may now compute the determinant (\ref{e:delta_multid}), i.e.,
the determinant of the restriction to ker $\cal Z$ of the linear map
$$
H_{\lambda, \xi}(\beta, \delta, \gamma) = 
\left(
\begin{array}{c} 
H^1_{\lambda, \xi}(\beta, \delta, \gamma) \\
H^2_{\lambda, \xi}(\beta, \delta, \gamma) \\
\end{array}
\right ) :=
\left (
\begin{array}{c}
\lambda(\bar \Omega \partial N + \bar N \partial \Omega )+ \sum_j i\xi_j \beta^1
\bar \Omega e_j  \\
\lambda \partial M + i \bar{\omega}\bar{X} \xi_1(\bar{M} \partial S + \partial Q- \sum_{j\ne1}
\delta^j \bar F^j) + 
\sum_{j \ne 1} i \bar{\omega}\bar{X}  \xi_j \partial F^j
\end{array}
\right ),
$$
which can be evaluated using an ingenious trick of \cite{Se.1} as
\ba \label{e:form}
\nn
\det (H_{\lambda, \xi} |_{\hbox {ker} \cal{Z}})&=&
C_1 \det  \left ( \begin{array}{c} H^1_{\lambda, \xi} \\
H^2_{\lambda, \xi} \\	\cal{Z}
 \end{array} \right )
=C_1C_2(\xi,\lambda) 
\det  \left ( \begin{array}{c} H^2_{\lambda, \xi} \\	\cal{Z}
 \end{array} \right )|_{\hbox {ker} H^1(\lambda, \xi)},\\
\ea
where $C_1= \det \Big( {\cal{Z}}|_{\hbox{ker} \cal{Z}^\perp}\Big) ^{-1}$
is well-defined thanks to full rank of $\cal{Z}$, assumption (H2)$'$,
and independent of $(\xi, \lambda)$, by the corresponding property of
$\cal{Z}$, but $C_2(\xi,\lambda)=\det H^1|_{\hbox{ker} (H^1)^\perp}$
and 
$$
\det  \left ( \begin{array}{c} H^2_{\lambda, \xi} \\	\cal{Z}
 \end{array} \right )|_{\hbox {ker} H^1(\lambda, \xi)}
$$
both depend on the specific dependence on $(\xi,\lambda)$
of the basis chosen for $\hbox {ker} H^1$.

Note that determinant (\ref{e:delta_multid}) is in the first place
defined only up to a constant factor depending on the parametrization
of $\setP$, so that we need only take care of the $(\xi,\lambda)$
dependence of $C_2$.
Rewriting 
\be \label{explicitH1}
H^1_{\lambda, \xi}(\beta, \delta, \gamma) =
\lambda(\sum_2^d \bar \omega \delta^j -\bar \omega^2 \gamma e_1)
+\sum_1^d i\xi_j \beta^1 \bomega e_j=0
\ee
as
\be
\label{e:rel1}
i\xi_1 \bar\omega \beta^1= \lambda \bar\omega^2\gamma,
\ee
and
\be
\label{e:rel2}
\lambda \bar\omega \delta^j=
-i\xi_j \bar\omega \beta^1, \quad j\ne 1,
\ee
and setting 
 as in \cite{Se.1}\footnote{We make the inessential change $\rho\to -\rho$
for convenience in later calculations.}
\be\label{e:comp1}
\beta^1=-\lambda \rho, \quad
\gamma=-i\bX \xi_1 \rho,
\ee
$\rho\in \C$, giving also the (new, multi-dimensional) relations
\be\label{e:comp2}
\delta^j= i\xi_j \rho,
\ee
and leaving $\beta^\ell$ free for $\ell\ne 1$, determines
a choice of basis for $\hbox{ker} H^1$, for which
$C_2(\xi,\lambda)$ has the simple form $\lambda^{d-1}$.

This fact is most easily verified by right-multiplying
$$
\det  \left ( \begin{array}{c} H^1_{\lambda, \xi} \\
H^2_{\lambda, \xi} \\	\cal{Z}
 \end{array} \right )
$$
by the determinant one matrix
$$
\left ( \begin{array}{cccc} 
0 & \alpha &\lambda\rho & 0\\
0_{2n-1} & 0_{2n-1} &0_{2n-1} & I_{2n-1}\\
I_{d-1} & 0_{d-1} &  i\tilde \xi & 0_{d-1}\\
0 & \beta & i\bar X \xi_1 \rho & 0\\
 \end{array} \right ),
$$
$\alpha i\bar X \xi_1\rho - \beta \lambda \rho=1$, to obtain
$$
\left ( \begin{array}{cc} 
N_1 &  0_{2n} \\
\hbox{*} &  N_2\\
 \end{array} \right ), \qquad 
N_1:=
\left ( \begin{array}{cccc} 
1 &&&  0 \\
0 &&& \omega \lambda I_{d-1}\\
 \end{array} \right ) 
$$
where $N_2$ is the $2n\times 2n$ matrix corresponding to
linear operator 
$
\left ( \begin{array}{cc} 
H^2\\
\cal{Z}\\
 \end{array} \right ) 
$
operating on $(\rho, \beta^2, \ldots, \beta^{2n})$ through
the compositions (\ref{e:comp1}) and (\ref{e:comp2}), and thus
\be
\hat \Delta(\xi, \lambda)= C_1\lambda^{d-1} \det N_2.
\ee

Alternatively, we may observe that $H^1$ is full rank whenever
$\lambda\ne 0$.  Observing also a posteriori that $\det N_2$ is
homogeneous degree $n+1$, we may conclude that
$C(\xi,\lambda)=\det H^1|_{\hbox{ker} (H^1)^\perp}$ as the ratio
of $n+d$ and $n+1$ degree homogeneous polynomials 
must be a constant times $\lambda^{d-1}$.
This discussion repairs a minor omission in \cite{Se.1},
where the dependence of $C_2$ on $(\xi, \lambda)$ is not explicitly discussed.

It remains to compute, under the compositions (\ref{e:comp1}), (\ref{e:comp2}),
the $2n\times 2n$ determinant $\det N_2$, which, transposing first and
second $n$-row blocks, may be expressed as
\be
\det N_2
=
\det \left (
\begin{array}{c}
\beta^1[\pd u/ \pd s] +\sum_2^{2n} \beta^l[w^l] + \gamma \bu' (0) +
\sum_2^d \delta^j [\pd u/\pd \delta^j]	\\   
\lambda \partial M + i\bar{\omega}\bar{X} \xi_1(\bar{M} \partial S + \partial Q- \sum_{j\ne1}
\delta^j \bar F^j) + 
\sum_{j \ne 1} i \bar{\omega}\bar{X}  \xi_j \partial F^j
\end{array}
\right ). 
\label{e:map}
\ee
Substituting from (\ref{e:comp1})--(\ref{e:comp2}) and the 
variational formulae (\ref{e:diff1})--(\ref{F1varsimple}),
and expressing $N_2$ as a matrix taking
$(\rho, \beta^2, \ldots, \beta^{2n})\to \C^{2n}$,
we obtain, similarly as in \cite{Se.1}, 
that $\det N_2=\hat \Delta \lambda^{1-d}$ 
is $\bar \omega^n$ times the determinant $\Delta_1$
defined in (\ref{e:leading}), giving
the desired relation $\Delta_1=\bar \omega^{-n} 
\lambda^{1-d}\hat \Delta$, and completing
the proof.

Namely, the first line of (\ref{e:map}) becomes
$$
-\lambda \rho [\pd u/ \pd s] + \sum_2^{2n} \beta^l[w^l] - 
i \bX \xi_1 \rho \bu' (0) + \sum_2^d i \xi_j \rho [\pd u/\pd \delta^j].
$$
The second line of (\ref{e:map}) becomes  
\ba
&& \lambda \bomega \left ( \gamma\bu(0)+\beta^1\int_0^{\bar{X}}
\frac{\partial u}{\partial s} dx_1+
\sum_2^{2n}\beta^l \int_0^{\bar{X}}w^l dx_1 +
\sum_2^d \delta^j
\int_0^{\bar{X}} \frac{\partial u}{\partial \delta^j} dx_1 
-\bar{M} \gamma \right ) 	\\	\nn
&+&  i \bomega\bX \xi_1 \left ( \bar{M} \beta^1 - \sum_2^{2n}
\beta^l \hL w^l - 
\sum_2^d  \delta^j F^j(\bu)  \right ) \\	\nn
&+& \sum_{j \ne 1} i \bomega\bX \xi_j \left ( \bomega \gamma f^j(\bu)(0)+
\bomega \int_0^{\bX} \pd f^j(\bu) d x_1 - \bomega \gamma F^j(\bu) \right ) \\ \nn
&=&
\lambda \bomega \left ( -i \bX \xi_1 \rho \bu(0)-\lambda \rho \int_0^{\bar{X}}
\frac{\partial u}{\partial s} dx_1+
\sum_2^{2n}\beta^l \int_0^{\bar{X}}w^l dx_1 +
\sum_2^d i \xi_j \rho
\int_0^{\bar{X}} \frac{\partial u}{\partial \delta^j} dx_1 
+\bar{M} i \bX \xi_1 \rho \right ) 	\\	\nn
&+&  i \bomega\bX \xi_1 \left ( -\bar{M} \lambda \rho - \sum_2^{2n}
\beta^l \hL w^l - \sum_2^d  i \xi_j \rho F^j(\bu)  \right ) \\	\nn
&+& \sum_{j \ne 1} i \bomega\bX \xi_j 
\left ( -\bomega i \bX \xi_1\rho f^j(\bu)(0)+
\bomega \int_0^{\bX} \pd f^j(\bu) d x_1 + 
\bomega i \bX \xi_1 \rho F^j(\bu) \right ) \\ \nn
&=&
-\bomega \rho \left (\lambda i \bX \xi_1 \bu(0) + \lambda^2 \int_0^{\bar{X}}
\frac{\partial u}{\partial s} dx_1 - \lambda i \sum_2^d \xi_j \int_0^{\bar{X}} 
\frac{\partial u}{\partial \delta^j} dx_1 
%
- \bX \xi_1 \sum_2^d \xi_j f^j(\bu)(0)  \right )  \\ \nn
&+& \bomega \sum_2^{2n} \beta^l \left ( \lambda \int_0^{\bX} w^l dx_1 
-i \bX \xi_1 \hL w^l \right ) + \bomega i \sum_2^d \xi_j \int_0^{\bX} \pd f^j(\bu) d x_1 \\ \nn
&=&
-\bomega \rho \left (- \lambda i \bX \xi_1 
\hL(\frac{\partial u}{\partial s})(0)+
\lambda^2 \int_0^{\bar{X}}
\frac{\partial u}{\partial s} dx_1 - \lambda i \sum_2^d \xi_j \int_0^{\bar{X}} 
\frac{\partial u}{\partial \delta^j} dx_1 - \bX \xi_1 \sum_2^d \xi_j \hL 
(\frac{\partial u}{\partial \delta^j})(0) \right ) \\ 	\nn
&+& \bomega \sum_2^{2n} \beta^l \left ( \lambda \int_0^{\bX} w^l dx_1 
-i \bX \xi_1 \hL w^l \right ) 
\\ \nn
&+& \bomega i \sum_2^d \xi_j \left ( 
-\lambda \rho \int_0^{\bX} 
\frac{\pd f^j}{\pd u} \frac{\pd u}{\pd s} dx_1 + \sum_2^{2n} \beta^l \int_0^{\bX} 
A^j w^l dx_1 + \sum_{k=2}^{d} i \xi_k \rho  \int_0^{\bX} 
\frac{\pd f^j}{\pd u} \frac{\pd u}{\pd \delta^k} dx_1 \right ) 
\ea
with (\ref{e:Mvar}) and other identities.
Denoting by ${\cal N}$ the matrix in (\ref{e:leading}) for which 
$\Delta_1=\det {\cal N}$, we find, 
comparing term by term, that the first $n$ rows of $N_2$
are equal to the first $n$ rows of ${\cal N}$,
while the last $n$ rows of $N_2$ are equal to $\bar \omega$ times
the last $n$ rows of ${\cal N}$.
Thus, $\det N_2= \bar \omega^n \det {\cal N}=\bar \omega^n \Delta_1$ 
as claimed, and we are done.


\begin{acknowledgment}
K.Z. thanks B. Texier for his interest in the problem, and
for several stimulating conversations.
Research of the authors was supported in part by the
National Science Foundation under Grants No. DMS-0204072 (M.O.)
and DMS-0300487 (K.Z.).
\end{acknowledgment}

\bibliographystyle{plain}

\end{document}